\documentclass[english, 12pt,leqno]{amsart}

\usepackage{a4,amsbsy,amsfonts,amsgen,
amsmath,
amsopn,amssymb,amsthm,amstext,enumerate,enumitem,exscale,
stmaryrd,epsfig,a4,graphicx,mathrsfs,relsize,url}
\usepackage{color}

\theoremstyle{definition}

\newcommand{\N}{\mathbf N}
\newcommand{\Z}{\mathbf Z}

\newcommand{\R}{\mathbf R}
\newcommand{\C}{\mathbf C}

\title
{On the prehistory of growth of groups}

\date{25 May 2022}


\author{Pierre de la Harpe}

\address{Pierre de la Harpe:
Section de math\'ematiques, 
Universit\'e de Gen\`eve,
C.P.~64, 
CH--1211 Gen\`eve 4.}
\email{Pierre.delaHarpe@unige.ch}

\begin{document}

\begin{abstract}
The subject of growth of groups has been active in the former Soviet Union since the early 50's
and in the West since 1968, when articles of \v{S}varc and Milnor have been published, independently.
The purpose of this note is to quote a few articles showing that,
before 1968 and at least retrospectively,
growth of groups has already played some role in various subjects.
\end{abstract}

\maketitle

The notion of growth for finitely generated groups appears in articles published independently
by Efremovich and \v{S}varc in the early 50's, and by Milnor in 1968
\cite{Efre--53, Svar--55, Miln--68a, Miln--68b}.
(\v{S}varc left Soviet Union in 1989, and now his name is rather written Schwarz.)
Very soon after his first paper, Milnor in \cite{Miln--68c} called attention to the fact that
\cite{Svar--55} ``contains many ideas utilized in [3]'' (where [3] = \cite{Miln--68b}).
\par

Before 1968, the paper \cite{Svar--55}, written by \v{S}varc during his undergraduate years,
was essentially ignored outside the former Soviet Union.
Concerning the period between 1955 and 1968,
we quote two extracts, from \cite[Definition I.6]{Avez--76}
and \cite[Item~0.5]{Grom--93}.
After having defined exponential growth and non-exponential growth, Avez writes:
``This notion is due to
V.\ Arnold (oral communication, 1965), \v{S}varc \cite{Svar--55}, and Milnor \cite{Miln--68b}.
Finite extensions of finitely generated nilpotent groups are the only known examples
of groups of non-exponential growth \cite{Wolf--68}.''
In what he calls some random historical remarks, Gromov writes:
``The ideas of the growth of balls, Folner sets and sets of conjugacy classes in groups
(especially in fundamental groups of manifolds of negative curvature,
see \cite{Marg--67} \cite{Marg--69})
were quite popular in the sixties among ergodic theorists in Moskow and Leningrad.
(Much of these ideas I learned at the time
from A.\ Vershik, D.\ Kazhdan and G.\ Margulis.)
Then the geometers took a part in the story and related the growth to curvature.
The first results here for non-negative curvature are due to A.\ \v{S}varc \cite{Svar--55}.
Similar results were obtained independently by J.\ Milnor \cite{Miln--68b}.''

\vskip.2cm

For a description of the results of Efremovich and \v{S}varc,
we quote the following lines from \cite{Svar--08}.
"My first serious work was inspired by Efremovich's remark that the
`volume invariant' of universal covering of a compact manifold
is a topological invariant of the manifold.
(If two compact manifolds are homeomorphic,
then the natural homeomorphism between universal coverings is uniformly continuous.
Efremovich proved that under certain conditions
the growth of the volume of a ball with radius tending to infinity
is an invariant of uniformly continuous homeomorphisms.)
I proved that the volume invariant of universal covering can be expressed
in terms of the fundamental group of the original manifold;
in modern language it is determined by the growth of the fundamental group.
I also gave estimates for volume invariants of manifolds with non-positive and with negative curvature. Thirteen years later J. Milnor published a paper containing the same results
with the only difference that Milnor was able to use in his proofs
some theorems derived after the appearance of my paper.
At the moment of writing his first paper in this direction Milnor did not know about my work,
but his second paper contained corresponding references.
The notion of growth of a group (volume invariant of a group in my terminology)
was studied later in numerous papers
(one should mention, in particular, the results by Gromov and Grigorchuk).
A new interesting field --- geometric group theory --- was born from these papers."
\par
For a short description of the work of Efremovich, see also \cite{Efremovich}.

\vskip.2cm

The importance of the subject of group growth was largely recognized with the results
of Gromov, showing that a finitely generated group has polynomial growth
if and only if it has a nilpotent subgroup of finite index \cite{Grom--81},
and Grigorchuk, showing the existence of groups of intermediate growth \cite{Grig--83, Grig--84}
(but it is still unknown whether there exist any finitely presented group of intermediate growth).
Several reviews of the subject have appeared,
of which we mention \cite{GrHa--97} and \cite{Grig--14},
and there is a nice exposition of the theory in the book by Mann \cite{Mann--12}.
Growth of groups extends naturally to the setting of locally compact groups;
in particular Guivarc'h and Jenkins \cite{Guiv--73, Jenk--73}
have characterized connected Lie groups of polynomial growth as those of Type (R),
i.e., as those for which ${\rm ad}(x)$ has purely imaginary eigenvalues
for all $x$ in the Lie algebra of the group
(this is considerably easier to prove than Gromov's characterization
of finitely generated groups of polynomial growth).
More generally, Losert has characterized locally compact groups (not necessarily Lie groups)
of polynomial growth in a series of papers, the first being \cite{Lose--87}.

\vskip.2cm

The purpose of this note is to mention a few articles
published before 1968, and for some even before 1955.
It can be seen retrospectively how the notion of group growth
has been used early for various purposes. 

\subsection{Carl Friedrich Gauss and the growth of $\Z^2$ (1834)}

The free abelian group of rank two, $\Z^2$, has to be seen
as the lattice of integer points in the Euclidean plane;
this has been so even before the concept of group was made precise in its present form.
Consider the length function on $\Z^2$ given by the Euclidean norm,
and the growth of $\Z^2$ as the function $R$ defined by
$$
R(t) = \left\vert \{ (a,b) \in \Z^2 \mid a^2 + b^2 \le t \} \right\vert
\hskip.5cm \text{for all} \hskip.2cm
t \ge 0 ,
$$
i.e., $R(t)$ is the number of points of $\Z^2$ in the disc of radius $\sqrt t$
centred at the origin.
The function $R(t)$ is interesting in number theory,
more precisely in the study of integers which are sums of two squares;
but we like to view also $R(t)$ as a function describing the growth of $\Z^2$.
In 1843, Gauss showed that
$$
\vert R(t) - \pi t \vert \le 2 \pi (1 + \sqrt{2t}) = O(\sqrt t) .
$$
In \cite[Pages 271 and 280]{Gauss}, Gauss wrote thirty values of $R(k)$,
including $R(10000) = 31417$ and $R(100000) = 314197$.
After Gauss, it has been shown that $\vert R(t) - \pi t \vert = O( t^\alpha )$
for values $\alpha < 1/2$, in particular for $\alpha = 1/3$ (Sierpinski, 1906);
the best estimate today seems to be $\vert R(t) - \pi t \vert = O( t^{\alpha + \varepsilon} )$
for $\alpha = 517/1648 = 0.31371...$ and for all $\varepsilon > 0$;
see \cite{BoWa}, as well as \cite{BeKZ--18}.
It is conjectured that the estimate holds for every $\alpha > 1/4$.
\par

For $k$ a non-negative integer, set
$r_2(k) = \left\vert \left\{ (a,b) \in \Z^2 \mid \sqrt{a^2+b^2} = k \right\} \right\vert$,
so that $R(k) = \sum_{j=0}^k r_2(j)$.
Values of $r_2(k)$ and $R(k)$ for small $k$ and relevant references
are given in \cite[A004018 and A057655]{OEIS}.
The series $\sum_{k=0}^\infty r_2(k) z^k$ is $(\theta_3(z))^2$,
where $\theta_3$ is the third Jacobi theta function
\cite[Chapter IV, Section 5]{CoSl--99}.
\par

\subsection{Word length and growth type of a finitely generated group}
\label{SectionGrowthType}

Let $\Gamma$ be a finitely generated group
and $S$ a finite generating set of $\Gamma$.
The \textbf{word length} function $\ell_S : \Gamma \to \N$ 
is defined by $\ell_S(\gamma) = \min \{ k \ge 0 \mid \gamma \in (S \cup S^{-1})^k \}$.
Let $\sigma(\Gamma, S; k)$ denote the cardinal
of the sphere $\{ \gamma \in \Gamma \mid \ell_S(\gamma) = k \}$
and $\beta(\Gamma, S; k)$ denote the cardinal
of the ball $\{ \gamma \in \Gamma \mid \ell_S(\gamma) \le k \}$.
It is straightforward to check that
$\sigma(\Gamma, S; k) \le
\vert S \cup S^{-1} \vert (\vert S \cup S^{-1} \vert - 1)^{k-1}$ for all $k \ge 1$;
in particular, $\beta(\Gamma, S; k) \le e^{bk}$
for an appropriate constant $b > 0$ and for all $k \ge 0$.
The group $\Gamma$ is said to be of 
\textbf{exponential growth} if there exist constants $c > 1$, $R \ge 0$
such that $\beta(\Gamma, S; k) \ge e^{ck}$ for all $k \ge R$,
of \textbf{subexponential growth} otherwise,
of \textbf{polynomial growth}
if there exist constants $C > 0$ and $d \in \N$ such that
$\beta(\Gamma, S; k) \le C k^d$ for all $k \ge 0$,
and of \textbf{intermediate growth} if it is of subexponential growth
and not of polynomial growth.
The definitions do not depend on the choice of $S$,
because the inequalities hold for one finite generating set $S$ if and only if
they hold for all finite generating sets.
\par

Word lengths, spheres and balls can be found in the literature
much before the theory of group growth.
For example $\ell_S(\gamma)$ appears
as the ``exponent of the substitution $\gamma$'' in \cite[Page 11]{Poin--82},
the paper in which Poincar\'e shows a presentation of a Fuchsian group
in terms of one of its fundamental polygons in the hyperbolic plane.
The word metric on $\Gamma$,
defined by $d_S(\gamma, \gamma') = \ell_S(\gamma^{-1}\gamma)$,
has been used systematically by Dehn in his first paper on decision problems in group theory;
see \cite{Dehn--11} and \cite[Pages 130 and 143]{DeSt--87}.
Spheres an balls, noted respectively $\Gamma_k$ and $\bigcup_{j=0}^k \Gamma_k$
appear in \cite{ArKr--63}, where the authors establish equidistribution in the $2$-sphere
of the points of the orbit of a semigroup generated by two appropriate rotations.

\subsection{Waclaw Sierpinski (1946), Georgii Adel'son--Vel'skii and
Yuli Anatoljevitch Shreider (1957),
Joseph Rosenblatt (1974), 
 and the supramenability of groups of subexponential growth}

In the 1929 paper which created the subject of amenability \cite{vNeu--29},
John von Neumann considers actions of a group $\Gamma$
on a set $X$ given with a non-empty subset $E$.
Such an action is \textbf{amenable} if there exists a finitely additive positive measure $\mu$ on $X$
normalized by $\mu(E) = 1$ and invariant by $\Gamma$
(the measure need not be finite, except of course when $E = X$).
The group $\Gamma$ itself is \textbf{amenable}
(eine messbare Gruppe in \cite{vNeu--29})
if every action of $\Gamma$
on every set $X$ given with $E = X$ is amenable,
and this holds as soon as the left action of $\Gamma$ on itself is amenable
(with $E = X = \Gamma$).
The group $\Gamma$ is \textbf{supramenable} (a~terminology due to Rosenblatt \cite{Rose--74}) 
if every action of $\Gamma$
on a set $X$ given with any subset $E \ne \emptyset$ is amenable,
and this holds as soon as the left action of $\Gamma$ on itself, with any $E$, is amenable.
The $\Gamma$-set $E$ has a \textbf{paradoxical decomposition} if there exist
a partition of $E$ in disjoint sets $A_1, \hdots, A_k, B_1, \hdots, B_\ell$
and elements $g_1, \hdots, g_k, h_1, \hdots, h_\ell$ in $\Gamma$
such that $E$ is equal to both the disjoint unions
$\bigsqcup_{i=1}^k g_i A_i$ and $\bigsqcup_{j=1}^\ell h_j B_j$.
A paradoxical decomposition of $E$
is an obstruction to the existence of $\mu$ as above \cite[Page 82]{vNeu--29};
remarkably it is the only obstruction: either $E$ has a paradoxical decomposition
or there exists a $\Gamma$-invariant finitely additive positive measure $\mu$ on $X$
normalized by $\mu(E) = 1$ \cite{Tars--36}.
\par

For example, the action on $X = \R^d$ of the isometry group of the Euclidean space $\R^d$
given with the unit ball $E$
is amenable when $d=1$ and $d= 2$, and is not when $d \ge 3$.
In dimension $3$, Hausdorff and Banach \& Tarski have obtained
famous results which express non-amenability in a spectacular way:
the action of the rotation group ${\rm SO}(3)$ on the unit ball in $\R^3$
is non-amenable, see \cite[Appendix to Chapter X, Page 469]{Haus--14},
and two bounded subsets $A$ and $B$ of $\R^3$
with non-empty interiors are equidecomposable
(this means that there exist partitions
$A = \bigsqcup_{i=1}^k A_i$, $B = \bigsqcup_{i=1}^k B_i$,
and isometries $g_1, \hdots, g_k$ of $\R^k$
such that $g_1A_1 = B_1, \hdots, g_kA_k = B_k$), see \cite{BaTa--24}.
\par

In \cite{Sier--46}, Sierpinski saw that any finitely generated subgroup
of the isometry group of $\R$ is of subexponential growth (indeed of polynomial growth),
and that this implies that the action of this isometry group on $\R$
is not paradoxical.
The argument shows essentially that the isometry group of $\R$ is supramenable,
and much more (see below).
\par

In \cite[Theorem 2]{AdSr--57},
it is shown that a finitely generated group of subexponential growth
is amenable.
\par

Later, Rosenblatt showed much more.
He introduced the terminology ``supramenable'', as defined above;
moreover he defined a group to be \textbf{exponentially bounded}
if all its finitely generated subgroups are of subexponential growth.
He showed that exponentially bounded groups are supramenable.
Moreover a finitely generated solvable group
either has a nilpotent subgroup of finite index,
and thus is of polynomial growth and supramenable,
or is not supramenable and contains a free semigroup on two generators,
and thus is of exponential growth \cite{Rose--74}.
\par

It is unknown whether there exist finitely generated groups of exponential growth
which are supramenable.

\vskip.2cm

We reproduce now Sierpinski's argument,
cast in the more general situation of a group $\Gamma$ acting on a set $X$
(instead of the affine group of $\R$ acting on $\R$),
and a nonempty subset $E$ of $X$.
Suppose that there exists a paradoxical decomposition of $E$:
there exist as above
subsets $A_1, \hdots, A_k, B_1, \hdots, B_\ell$ of $E$
and a subset $S$ of elements $g_1, \hdots, g_k, h_1, \hdots, h_\ell$ of~$\Gamma$
(not necessariliy distinct from each other) such that
$$
E = \bigg( \bigsqcup_{i=1}^k A_i \bigg) \sqcup \bigg( \bigsqcup_{j=1}^\ell B_i \bigg)
= \bigsqcup_{i=1}^k g_i A_i = \bigsqcup_{j=1}^\ell h_j B_j .
$$
The following argument shows that the subgroup of $\Gamma$ generated by $S$
has exponential growth.
\par

Set $A = \bigsqcup_{i=1}^k A_i$, $B = \bigsqcup_{j=1}^\ell B_j$.
Define bijections $\varphi \, \colon E \to A$ and $\psi \, \colon E \to B$
by $\varphi (x) = g_i^{-1}(x)$ when $x \in g_i A_i$
and $\psi (x) = h_j^{-1} (x)$ when $x \in h_j B_j$.
Choose $x_0 \in E$.
Observe first that $\varphi (x_0) \ne \psi (x_0)$, because $A$ and $B$ are disjoint,
then that $\varphi \varphi (x_0), \varphi \psi (x_0), \psi \varphi (x_0), \psi \psi (x_0)$
are also distinct, because $\varphi$ and $\psi$ are injective and $A \cap B = \emptyset$,
and so on.
This shows that, for any positive integer $k$, the $2^k$ words of length $k$ in $\varphi$ and $\psi$
are maps $E \to E$ with distinct values at $x_0$.
For any of these words, say $\chi$, the value $\chi (x_0)$
is of the form $s_1^{-1} s_2^{-1} \cdots s_k^{-1} (x_0)$,
for $s_1, s_2, \hdots, s_k \in S$.
It follows that the subgroup of $\Gamma$ generated by $S$
has at least $2^k$ distinct elements $\gamma$ of word length $\ell_S (\gamma) \le k$,
and this ends the argument.
\par

Sierpinski's argument shows that a group $\Gamma$ which can act on a pair $X \supset E$
such that $E$ has a paradoxical decomposition
has a finitely generated subgroup of exponential growth.
By contraposition, it follows that
a finitely generated group of subexponential growth is supramenable,
a result much stronger than the one in \cite{AdSr--57},
and a proof much more direct than the one in \cite{Rose--74}.

\subsection{Hans Ulrich Krause
and finitely generated abelian groups with isomorphic Cayley graphs (1953)}

In his thesis \cite[Satz 16.1]{Krau--53},
Krause shows that two finitely generated abelian groups have isomorphic Cayley graphs
with respect to well-chosen generating sets
if and only if the two following conditions are satisfied:
(i) the two groups have the same rank, and
(ii) their torsion groups have the same order.
In the proof, it is shown that the rank of a finitely generated abelian group
$\Gamma$ is the polynomial growth rate $\lim_{k \to \infty} (\ln \vert S^k \vert ) / (\ln k)$,
where $S$ is a symmetric generating set of $\Gamma$.

\subsection{Jacques Dixmier
and polynomial growth of nilpotent connected Lie groups (1960, 1966)}

Lemma 3 of \cite{Dixm--60} is the following.
Let $G$ be a nilpotent connected Lie group, $\mu$ a Haar measure on $G$,
and $H$ a compact subset of $G$.
Then there exists an integer $N$ (which depends on $G$ but not on $H$)
such that $\mu(H^k) = O(k^N)$ when $k \to \infty$;
in other words, $G$ is a group of polynomial growth.
\par

The lemma is used by Dixmier in the proof of the following result.
Consider a locally compact group $G$,
the group algebra $L^1(G)$,
and the two-sided ideal $I$ of those elements $f \in L^1(G)$ such that,
for every irreducible unitary representation $\pi$ of $G$,
the operator $\pi(f)$ is of finite rank.
If $G$ is a nilpotent connected Lie group, then $I$ is dense in $L^1(G)$.
(The same property of $I$ was established earlier for semisimple Lie groups by Harish--Chandra.)

Polynomial growth has been established later for
solvable connected Lie groups of type (R), in \cite{Dixm--66}.

\subsection{Henri Dye and orbital equivalence (1963)}

Theorem 1 of \cite{Dye--63} establishes the following.
Let $\Gamma$ be a finitely generated group, generated by a finite subset $F$.
The notation of Dye is
$h_1 = \vert F \vert$ and $h_k = \vert F^k \smallsetminus F^{k-1} \vert$ for $k \ge 2$.
If
$$
\inf_{k \ge 1} \frac{h_{2k}}{h_1 + \cdots + h_k} = 0 ,
$$
then $\Gamma$ is approximately finite.
In particular, finitely generated groups of polynomial growth are approximately finite.
\par

To define approximate finiteness,
consider actions of countable groups on non-atomic standard probability spaces
by measure preserving transformations.
Two such actions of $\Gamma_1$ on $X_1$ and $\Gamma_2$ on $X_2$
are orbit equivalent if there exists a measure preserving Borel isomorphism
$f \ \colon X_1 \to X_2$ such that $f(\Gamma_1 x)$ coincides with the orbit $\Gamma_2 f(x)$
for almost all $x$ in $X_1$.
Consider some ergodic measure preserving action of the infinite cyclic group $\Z$
on a non-atomic standard probability space;
a basic example is the Bernoulli shift action $\beta$ of $\Z$ on $(\Z / 2\Z)^\Z$.
A countable group $\Gamma$ is approximately finite in the sense of Dye if,
for every ergodic measure preserving action $\alpha$ of $\Gamma$
on a non-atomic probability space $X$,
the actions $\alpha$ and $\beta$ are orbit equivalent.
\par

It is now known that an infinite countable group is approximately finite
if and only if it is amenable \cite{OrWe--80, Hjor--05}.

\subsection{Gregori Margulis, growth of fundamental group
and existence of Anosov flows (1967)}

On a compact Riemannian smooth manifold $M$,
an Anosov flow is a smooth flow $\Phi = \{ \Phi_t \}_{t \in \R}$
which satisfies the following conditions.
There exists a $\Phi$-invariant continuous splitting
$TM = E^u \oplus E^T \oplus E^s$
of the tangent bundle of $M$,
where the three terms are respectively
the unstable (or expanding) subbundle of $TM$,
the line bundle tangent to $\Phi$,
and the stable (or contracting) subbundle of $TM$,
and there exist constants $\nu > 0$, $c > 0$, such that
\par

$\Vert ({\Phi_t})_* (v) \Vert \ge c e^{\nu t} \Vert v \Vert$
and
$\Vert ({\Phi_{-t}})_* (v) \Vert \le c e^{- \nu t} \Vert v \Vert$
for all $v \in E^u$ and $t \ge 0$
\par

$\Vert (\Phi_t)_* (v) \Vert \le c e^{-\nu t} \Vert v \Vert$
and
$\Vert (\Phi_{-t})_* (v) \Vert \ge c e^{\nu t} \Vert v \Vert$
for all $v \in E^s$ and $t \ge 0$
\par
\noindent
(the two conditions with $({\Phi_{-t}})_* (v)$
follow from the two conditions with $({\Phi_t})_* (v)$,
see \cite[Page 121]{AnSi--67}).
\par

In one of his first published papers, Margulis shows that, if a $3$-dimensional manifold $M$
has an Anosov flow, then the fundamental group of $M$ has exponential growth
\cite{Marg--67}.
This has been generalized to manifolds of higher dimensions
and Anosov flows with one of the subbundles $E^u, E^s$
of rank one \cite{PlTh--72}.
\par

For the contrast, let us quote the following result of Franks.
On a compact Riemannian smooth manifold $M$,
a $\mathcal C^1$ map $f \, \colon M \to M$ is expanding
if there are constants $\lambda > 1$ and $c > 0$
such that $\Vert T f^m v \Vert \ge c \lambda^m \Vert v \Vert$
for all $v \in TM$ and $m \ge 1$.
Here is the result:
If a compact manifold admits an expanding self-map,
then its fundamental group has polynomial growth
\cite[Theorem~8.3]{Fran--70}.

\subsection{Harry Kesten and recurrent random walks on groups (1967)}

Let $\Gamma$ be a finitely generated group.
A symmetric probability measure $\mu$ on $\Gamma$
such that $\{ \gamma \in \Gamma \mid \mu(\gamma) > 0 \}$
is a finite generating set gives rise to a random walk on $\Gamma$.
The group $\Gamma$ is recurrent
if there exists such a measure such that the associated random walk is recurrent
(equivalently if for any such measure the associated random walk is recurrent).
It is a classical theorem of P\'olya that the simple random walk on $\Z^d$
is recurrent if $d \le 2$ and transient if $d \ge 3$ \cite{Poly--21},
and it has been known since at least 1962 that an infinite finitely generated abelian group
is recurrent if and only if it is a finite extension of $\Z$ or a finite extension of $\Z^2$
\cite{Dudl--62}.
\par

Kesten has conjectured that the recurrence of a group depends on its growth type.
It is conjectured more precisely in \cite[Conjecture 4]{Kest--67}
that a finitely generated group which is recurrent cannot be of exponential growth.
The conjecture was made more precise (the growth of a recurrent group is at most quadratic)
and generalized to second countable locally compact groups;
see the introduction of \cite{GuRa--12}.
For discrete groups, the final result is due to Varopoulos (1986) :
a finitely generated group is recurrent if and only if it is of at most quadratic growth,
if and only it is either finite, or a finite extension of $\Z$, or a finite extension of $\Z^2$;
see \cite{VaSC--92}.
(A finitely generated group $\Gamma$ is of quadratic growth if,
with the notation of Section~\ref{SectionGrowthType}, there exists a constant $C > 0$
such that $\beta(\Gamma, S; k) \le Ck^2$ for all $k \ge 0$.)

\subsection{Generating functions}
\label{SectionGenF}

To encode a sequence $(a_k)_{k \ge 0}$ of integral numbers,
several types of series or functions can be used,
and the best choice depends on the subject.
One choice is the
\textbf{ordinary generating function} of the sequence $(a_k)_{k \ge 0}$:
$$
\Sigma(z) = \sum_{k=0}^\infty a_k z^k \in \Z [[z]] .
$$
When $\Sigma(z)$ converges for $z$ small enough
and whenever possible,
we like to identify the ``simple function of analysis'' of which
$\Sigma(z)$ is the Taylor series at the origin.
The book \cite{FlSe--09} is a very rich source of examples and theorems on these generating functions.
\par

An early example occurs in
a letter of Euler to Goldbach dated September~4, 1751.
The letter is reproduced partly in \cite[Section~I.1]{FlSe--09},
and in full in \cite[Letter 154, Pages 489--491]{Euler}.
(In \cite[Section 1.2.9]{Knut--97},
Knuth mentions three earlier appearances of generating functions
by de Moivre, by Stirling, and by Euler in connection with
numbers of partitions of integers.)
In his letter, Euler considers the number $c_k$ of decompositions
of a convex $(k+2)$-gon in triangles; set moreover $c_0 = 1$.
The generating function of $(c_k)_{k \ge 0}$
$$
\begin{aligned}
\sum_{k \ge 0} c_k z^k &= 1 + z + 2z^2 + 5z^3 + 14z^4 + 42z^5 + 132 z^6 + 429 z^7 + \cdots
\\
&= 1 + \sum_{k=1}^\infty
\frac{ 2 \cdot 6 \cdot 10 \cdot \hdots \cdot (4k-2) }{2 \cdot 3 \cdot 4 \cdot \hdots \cdot \hskip.2cm (k+1)}
z^k
= \frac{ 1 - \sqrt{1-4z} }{ 2z } 
\end{aligned}
$$
is algebraic.
(Euler does not consider our $c_0$, and writes $a$ instead of our $z$,
so that his series sums up to
$\frac{ 1 - 2a - \sqrt{1-4a} }{ 2aa }$.)
The numbers $c_k$ are now known as Catalan numbers,
and are often written in terms of binomial coefficients: $c_k = \frac{1}{k+1} \binom{2k}{k}$.
For $214$ different kinds of objects that are counted using Catalan numbers
and for a historical survey, see \cite{Stan--15}.
\par

The simplest sequences
are those which satisfy a linear recurrence relation;
they correspond precisely to rational generating functions.
More precisely, consider a positive integer $d$
and complex numbers $q_1, q_2, \hdots, q_d$ with $q_d \ne 0$.
Set $Q(z) = 1 + q_1 z + q_2 z^2 + \cdots + q_d z^d
= \prod_{j=1}^e (1 - \gamma_j z)^{d_j}$,
where $\gamma_1, \hdots, \gamma_e \in \C$ are the distinct roots of $Q$
and $d_1, \hdots, d_e$ their multiplicities;
note that $\sum_{j=1}^e d_j = d$.
Then, for a sequence $(a_k)_{k \ge 0}$,
the following conditions are equivalent
\begin{enumerate}[label=(\roman*)]
\item[(R1)]
$\sum_{k \ge 0} a_k z^k = P(z) / Q(z)$
for some polynomial $P(z)$ of degree less than~$d$,
\item[(R2)]
$a_{k+d} + q_1 a_{k+d-1} + q_2 a_{k+d-2} + \cdots + q_d a_k = 0$
for all $k \ge 0$,
\item[(R3)]
$a_k = \sum_{j=1}^e P_j(k) \gamma_j^k$
for all $k \ge 0$, for some polynomials $P_j (z)$ of degree less than $d_j$
(with $j = 1, \hdots, e$).
\end{enumerate}
For this, and for variations (when $\deg P \ge d$ or when $Q(z) = (1-z)^d$),
see \cite[Chapter 0]{Stan--96}.
\par

The Fibonacci sequence $(F_k)_{k \ge 0} = (0, 1, 1, 2, 3, 5, 8, 13, 21, 34, \hdots)$
is a notorious example:
\begin{enumerate}[label=(\roman*)]
\item[(1)]
generating function $\sum_{k=0}^\infty F_k z^k = \frac{z}{1 - z -z^2}$,
\item[(2)]
linear recursion $F_{k+2} - F_{k+1} - F_k = 0$ for all $k \ge 0$,
\item[(3)]
and Binet formula
$F_k = \frac{1}{ \sqrt 5} \left( \frac{1 + \sqrt 5 }{2} \right)^k
- \frac{1}{ \sqrt 5} \left( \frac{1 - \sqrt 5 }{2} \right)^k$
(published by Binet in 1845, but already in 
\cite[\S~7]{Ber--1728} and \cite[Page 128]{Eul--1767}).
\end{enumerate}

\subsection{Growth series for finitely generated groups}
\label{SectionGrowthSGps}

Let $\Gamma$ be a finitely generated group
and $S$ a finite generating set of $\Gamma$.
For each $k \ge 0$, let $\sigma(\Gamma, S; k)$ be the cardinal of the sphere of radius $k$
and let $\beta(\Gamma, S; k)$ the cardinal of the ball of radius $k$,
as defined in Section~\ref{SectionGrowthType}.
The \textbf{growth series} of the pair $(\Gamma, S)$ is the generating function
$$
\Sigma(\Gamma, S; z)
= \sum_{k=0}^\infty \sigma(\Gamma, S; k) z^k
= \sum_{\gamma \in \Gamma} z^{\ell_S(\gamma)}
\in \Z [[z]] .
$$
The radius of convergence of this series is strictly positive and is $1/\omega(G, S)$,
where $\omega(G, S) = \lim_{k \to \infty} \sigma(\Gamma, S; k)^{1/k}$
is the exponential growth rate of the pair $(G, S)$.
It is sometimes better to consider
$$
B(\Gamma, S; z) 
= \sum_{k=0}^\infty \beta(\Gamma, S; k) z^k
= \frac{\Sigma(\Gamma, S; z)}{1-z} .
$$
\par

For example,
for the infinite cyclic group $\Gamma = \Z$ generated by $S = \{1\}$, we have
$$
\boxed{
\Sigma(\Z, \{1\}; z) = 1 + 2z + 2z^2 + 2z^3 + 2z^4 + 2 z^5 + \cdots = \frac{1+z}{1-z} .
}
$$
More generally, for the free abelian group $\Z^n$ generated by a basis $S_n$,
we have
$$
\Sigma(\Z^n, S_n; z)
= \sum_{k=0}^\infty \Bigg( \sum_{\ell=0}^n \binom{n}{\ell} \binom{k+n-\ell-1}{k-\ell} \Bigg) z^k
= \left( \frac{1+z}{1-z} \right)^n .
$$
The infinite sum simplifies to
$\Sigma_{k=0}^\infty 4k z^k$,
$\Sigma_{k=1}^\infty (4k^2+2)z^i$,
$\Sigma_{k=1}^\infty \frac{8k(k^2+2)}{3} z^k$,
when $n = 2, 3, 4$, respectively
(sequences A008574, A005899, A008412 in \cite{OEIS}).
\par

Growth series have been studied for several other classes of groups.
For a Coxeter system $(\Gamma, S)$ with $S$ finite, 
the growth series $\Sigma(\Gamma, S; z)$ is a rational function.
See exercise 26 of Chap.\ IV \S~1 and exercise 10 of Chap.\ VI \S~4
in \cite{Bour--68}.
This function has interesting values; for example,
its value at $1$ is rational
and is the inverse of the Euler--Poincar\'e characteristic of the group $\Gamma$
\cite[Proposition 17, Page 112]{Serr--71}.

For a Gromov hyperbolic group $\Gamma$ and an arbitrary generating set $S$,
Gromov has shown that $\Sigma(\Gamma, S; z)$ is a rational function
\cite[Corollary 5.2.A']{Grom--87}.
This generalizes a result of Cannon \cite[Theorem~7]{Cann--84}.
\par

There are some groups $\Gamma$ with generating sets $S$
such that $\Sigma(\Gamma, S; z)$ is an irrational algebraic function
\cite{Parr--92}.
The growth series of a pair $(\Gamma, S)$ can also be transcendental.
Stoll showed that there are groups $\Gamma$ with two finite generating sets $S, T$
such that $\Sigma(\Gamma, S; z)$ is rational and $\Sigma(\Gamma, T; z)$ transcendental
\cite{Stol--96}.

\subsection{Hilbert series}
\label{SectionHilbertSeries}

Consider again the group $\Gamma = \Z^n$ for some $n \ge 1$
and an arbitrary finite generating set $S$.
Then there exists a polynomial $P \in \Z[z]$ such that
$$
\Sigma(\Z^n, S; z) = \frac{P(z)}{ (1-z)^n } .
$$
Here is one way to show this: the group algebra $\C [\Gamma]$,
with linear basis $(\delta_\gamma)_{\gamma \in \Gamma}$
and multiplication defined by $\delta_\gamma \delta_{\gamma'} = \delta_{\gamma \gamma'}$,
has a filtration $\C[\Gamma] = \bigcup_{k \ge 0} B_k$
where $B_k$ is the linear subspace generated by $\{ \delta_\gamma \mid \ell_S(\gamma) \le k \}$;
set moreover $B_{-1} = \{0\}$.
The associated graded algebra $A = \bigoplus_{k \ge 0} (B_k / B_{k-1})$
is commutative and generated by a finite set of elements of degree $1$.
It is a theorem of Hilbert that the Hilbert series
$$
\sum_{k \ge 0} \dim_\C (B_k / B_{k-1}) z^k = \Sigma(\Z, S; z)
$$
of such an algebra is rational of the form $\frac{P(z)}{ (1-z)^n }$;
for a proof, see for example \cite[Theorem 11.1]{AtMa--69}.
The observation that the growth series of $(\Gamma, S)$
is the Hilbert series of an appropriate graded algebra,
and thus in particular a rational function,
is due to several authors, including \cite{Wagr--82}.
\par

The ``theorem of Hilbert'' refers to Theorem IV in \cite[Page 512]{Hilb--90}.
In fact, Hilbert shows that the series
satisfies a condition like (R3) of our Section~\ref{SectionGenF}, rather than (R1).
But it was already standard in this time to write
``Hilbert series'' which are rational functions
for the dimensions of the homogeneous components of a graded algebra;
I am grateful to Hanspeter Kraft for showing me that this can be found
in the work of Sylvester on the theory of invariants, around 1880; 
see for example papers 38, 40, and 59, in \cite{Sylvester}.
\par

Hilbert series are also called Poincar\'e series, especially when
they encode dimensions of homology spaces; see \cite{Babe--86}.
\par

More generally, when $\Gamma$ is a virtually abelian finitely generated group
and $S$ an arbitrary finite generating set,
the series $\Sigma(\Gamma, S; z)$ is rational \cite{Bens--83}.

\subsection{Eug\`ene Ehrhart
and the number of integral points in the multiples of a polytope (1962)}

\label{SectionEhrhart}
Consider an Euclidean space $V$ of dimension $n$,
a lattice $\Gamma$ in $V$,
i.e., a subgroup of $V$ isomorphic to $\Z^n$ and generated by a basis of~$V$,
a polytope $P$ which is the convex hull of a finite subset of $\Gamma$, 
and for each non-negative integer $k$ the number $E_P(k)$ of points in $kP \cap \Gamma$.
In 1962, Ehrhart published a note on the numbers $E_P(k)$
and the series $\sum_{k=0}^\infty E_P(k) z^k$ \cite{Ehrh--62, Brio--95}.
For a polytope of non-empty interior,
this series is a growth series of the group $\Gamma \approx \Z^n$
for a particular choice of generating sets.
\par

Note that, 
for the lattice $\Z^n$ in $\R^n$ and the convex hull $P = {\rm Conv}(\pm e_1, \hdots, \pm e_n)$,
where $\{e_1, \hdots, e_n\}$ is the standard basis of $\R^n$,
we have, with the notation of Section~\ref{SectionGrowthSGps},
$$
\sum_{k=0}^\infty E_P(k) z^k
= B(\Z^n, \Z^n \cap P; z) 
= \frac{1}{1-z} \left( \frac{1+z}{1-z} \right)^n .
$$
Other cases are studied from this point of view in \cite{BaHV--99}.
For example, when $V = \{ (x_1, \hdots, x_{n+1}) \in \R^{n+1} \mid \sum_{j=1}^{n+1} x_i = 0 \}$,
$\Gamma = \Z^{n+1} \cap V \approx \Z^n$,
and $P$ is the convex hull of $\{ \pm(e_i - e_j) \mid 1 \le i < j \le n+1 \}$,
$$
\sum_{k=0}^\infty E_P(k) z^k
= B(\Gamma, \Gamma \cap P, z)
= \frac{1}{ (1-z)^{n+1} } \sum_{j=0}^n \binom{n}{k}^2 z^j 
= \frac{1}{1-z} P_n \left( \frac{1+z}{1-z} \right) ,
$$
where $P_n$ is the Legendre polynomial of degree $n$.

\subsection{Theta functions}

Consider a Euclidean vector space $V$ of dimension $n$,
with scalar product denoted by $\langle \cdot \hskip-.1cm \mid \hskip-.1cm \cdot \rangle$,
and a lattice $\Gamma$ in $V$.
For elements of $\Gamma$, consider no longer the word length as above,
but rather the norm $\Gamma \to \R_+, \hskip.2cm
x \mapsto \Vert x \Vert = \sqrt{ \langle x \hskip-.1cm \mid \hskip-.1cm x \rangle }$.
The \textbf{theta function} of $\Gamma$ is defined by
$$
\Theta_\Gamma(\tau) = \sum_{x \in \Gamma} e^{i \pi \tau \Vert x \Vert^2} ,
$$
so that $\Theta_{\Gamma}$ is a holomorphic function on the upper half-plane
$\{ \tau \in \C \mid \text{Im}(\tau) > 0 \}$.
When $\Gamma$ is an integral lattice, namely when
$\langle x \hskip-.1cm \mid \hskip-.1cm y \rangle \in \Z$ for all $x, y \in \Gamma$,
the theta series is alternatively viewed as a power series in
$q = e^{i \pi \tau}$~:
$$
\Theta_\Gamma(q) = \sum_{x \in \Gamma} q^{ \Vert x \Vert^2 }
= \sum_{r=0}^\infty
\vert \{ x \in \Gamma \mid \langle x \hskip-.1cm \mid \hskip-.1cm x \rangle = r \} \vert \hskip.1cm q^r .
$$
For example,
when $\Gamma = \Z$ is embedded the standard way in the real line $V = \R$,
the series is
$$
\boxed{
\Theta_\Z(q) = 1 + 2q + 2q^4 + 2q^9 + 2q^{16} + 2 q^{25} + \cdots = \theta_3(q) 
}
$$
where $\theta_3$ is as above the third Jacobi theta function.
More generally, for $\Z^n$ embedded the standard way
in the standard Euclidean space $\R^d$,
we have $\Theta_{\Z^n} (q) = \left( \theta_3 (q) \right)^n$
\cite[\emph{op.\ cit.}]{CoSl--99}.
\par

It is tempting to compare the two boxed formulas of this paper related to $\Z$,
and more generally to speculate whether theta functions
could be of some interest for other groups than lattices in Euclidean spaces.

\vskip.3cm

I am grateful to Jean-Paul Allouche, Ivan Babenko, Slava Grigorchuk, Avinoam Mann,
Hanspeter Kraft, Nicolas Monod, Tatiana Nagnibeda, and Igor Pak
for useful comments and suggestions.


\begin{thebibliography}{AkGW--08}

\bibitem[AdSr--57]{AdSr--57}
G.M.\ Adel'son--Vel'skii and Yu.A. Sreider,
\emph{The Banach mean on groups}.
Uspehi Mat.\ Nauk (N.S.) \textbf{12} (1957), no.\ 6(78), 131--136.

\bibitem[AnSi--67]{AnSi--67}
D.V.\ Anosov and Ya.G.\ Sinai,
\emph{Some smooth ergodic systems}.
Russian Math.\ Surveys \textbf{22} (1967), no.\ 5, 103--167.

\bibitem[ArKr--63]{ArKr--63}
V.I.\ Arnol'd and V.I.\ Krylov,
\emph{Uniform distribution of points on a sphere and certain ergodic properties of solutions
of linear ordinary differential equations in a complex domain}.
Dokl.\ Akad.\ Nauk SSSR \textbf{148} (1963), 9--12.

\bibitem[AtMa--69]{AtMa--69}
M.F.\ Atiyah and I.G.\ Macdonald,
\emph{Introduction to commutative algebra}.
Addison--Wesley, 1969.

\bibitem[Avez--76]{Avez--76}
A.\ Avez,
\emph{Croissance des groupes de type fini et fonctions harmoniques}.
In ``Th\'eorie ergodique (Actes journ\'ees ergodiques, Rennes, 1973/1974)'',
Springer Lecture Notes in Math.\ \textbf{532} (1976), 35--49.

\bibitem[Babe--86]{Babe--86}
I.K.\ Babenko,
\emph{Problems of growth and rationality in algebra and topology}.
Russian Math.\ Surveys \textbf{41} (1986), no.\ 2, 117--175.


\bibitem[BaHV--99]{BaHV--99}
R.\ Bacher, P.\ de la Harpe, and B.\ Venkov,
\emph{S\'eries de croissance et polyn\^omes d'Ehrhart associ\'es aux r\'eseaux de racines},
Ann.\ Inst.\ Fourier \textbf{49} (1999), no.\ 3, 727--762.


\bibitem[BaTa--24]{BaTa--24}
 S.\ Banach and A.\ Tarski,
 \emph{Sur la d\'ecomposition des ensembles de points en parties respectivement
congruentes}.
Fund.\ Math.\ \textbf{6} (1924), 244--277
[Banach, Oeuvres, Vol.\ I, 118--148 et 325--327].


\bibitem[Bens--83]{Bens--83}
M.\ Benson,
\emph{Growth series of finite extensions of $\Z^n$ are rational}.
Invent.\ Math.\ \textbf{73} (1983), no.\ 1, 251--269.

\bibitem[BeKZ--18]{BeKZ--18}
B.C.\ Berndt, S.\ Kim, and A.\ Zaharescu,
\emph{The circle problem of Gauss and the divisor problem of Dirichlet --- still unsolved}.
Amer.\ Math.\ Monthly \textbf{125} (2018), no.\ 2, 99--114.

\bibitem[Ber--1728]{Ber--1728}
D.\ Bernoulli,
\emph{Observationes de seriebus
quae formantur ex additione
vel subtractione quacunque
terminorum se mutuo consequentium}.
Commentarii academiae scientiarum imperialis petropolitanae, tomus III,
annum $\subset$I$\supset$ !$\supset\subset\subset$ XXVIII, 85--100.
\url{https://archive.org/details/commentariiacade03impe/page/84/mode/2up?view=theater}


\bibitem[Bour--68]{Bour--68}
N.\ Bourbaki,
\emph{Groupes et alg\`ebres de Lie, chapitres 4, 5 et 6}.
Hermann, 1968.

\bibitem[BoWa]{BoWa}
J.\ Bourgain and N.\ Watt,
\emph{Mean square of zeta function, circle problem, and divisor problem revisited}.
ArXiv:1709.04340v1, 13 Sep 2017.


\bibitem[Brio--95]{Brio--95}
M.\ Brion,
\emph{Points entiers dans les polytopes convexes}.
S\'eminaire Bourbaki \textbf{780}, 1993--1994,
Ast\'erisque \textbf{227} (1995), 145--169.

\bibitem[Cann--84]{Cann--84}
J.W.\ Cannon,
\emph{The combinatorial structure of cocompact discrete hyperbolic groups}.
Geom.\ Dedicata \textbf{16} (1984), no.\ 2, 123--148.

\bibitem[CoSl--99]{CoSl--99}
J.H.\ Conway and N.J.\ Sloane,
\emph{Sphere packings, lattices and groups}, Third edition.
Springer, 1999.

\bibitem[Dehn--11]{Dehn--11}
M.\ Dehn,
\emph{\"Uber unendliche diskontinuierliche Gruppen}.
Math.\ Ann.\ \textbf{71} (1912), 116--144
(Volume 71 is dated 1912, but the first part, Pages 1-144,
appeared on July 25, 1911).
See also
\emph{On infinite discontinuous groups},
fourth article in \cite{DeSt--87}.

\bibitem[DeSt--87]{DeSt--87}
M.\ Dehn,
\emph{Papers on group theory and topology,
translated and introduced by John Stillwell}.
Springer, 1987.

\bibitem[Dixm--60]{Dixm--60}
J.\ Dixmier,
\emph{Op\'erateurs de rang fini dans les repr\'esentations unitaires}.
Inst.\ Hautes \'Etudes Sci.\ Publ.\ Math.\ {\bf 6} (1960), 13--25.

\bibitem[Dixm--66]{Dixm--66}
J.\ Dixmier,
\emph{Sur les groupes de Lie r\'esolubles \`a racines purement imaginaires}.
Bull.\ Sci.\ Math.\ (2) \textbf{90} (1966), 5--16.

\bibitem[Dudl--62]{Dudl--62}
R.M.\ Dudley,
\emph{Random walks on abelian groups}.
Proc.\ Amer.\ Math.\ Soc.\ \textbf{13} (1962), 447--450.

\bibitem[Dye--63]{Dye--63}
H.A.\ Dye,
\emph{On groups of measure preserving transformations. II}.
Amer.\ J.\ Math.\ \textbf{85}, no.\ 4 (1963), 551--576.

\bibitem[Efre--53]{Efre--53} 
V.A.\ Efremovich
\emph{On proximity geometry of Riemannian manifolds}.
Uspekhi Math.\ Nauk.\ \textbf{8} (1953), 189--191. 
[English translation:
Amer.\ Math.\ Transl.\ (2) \textbf{39} (1964), 167--170.]


\bibitem[Efremovich]{Efremovich}
\emph{Vadim Arsen'evich Efremovich (obituary)},
by V.S.\ Makarov, A.A.\ Mal'tsev, S.P.\ Novikov, S.S.\ Ryshkov, A.S.\ Schwarz.
Russian Math.\ Surveys \textbf{45} (1990), no.\ 6, 137--138.

\bibitem[Ehrh--62]{Ehrh--62}
E.\ Ehrhart,
\emph{Sur les poly\`edres rationnels homoth\'etiques \`a $n$ dimensions}.
C.R.\ Acad.\ Sci.\ Paris \textbf{254} (1962), 616--618.


\bibitem[Eul--1767]{Eul--1767}
L.\ Eulero,
\emph{Observationes analyticae}.
The Euler Archive,
\hfill\par\noindent
\url{https://scholarlycommons.pacific.edu/euler-works/326/}

\bibitem[Euler]{Euler}
Leonhardi Euleri,
\emph{Opera omnia. Series 4 A. Commercium epistolicum. Vol. 4.1. Leonhardi Euleri commercium epistolicum cum Christiano Goldbach. Pars I / Correspondence of Leonhard Euler with Christian Goldbach. Part I.}
Original texts in Latin and German. Edited by F.\ Lemmermeyer and M.\ Mattm\"uller.
Springer, Basel, 2015.


\bibitem[FlSe--09]{FlSe--09}
P.\ Flajolet and R.\ Sedgewick,
\emph{Analytic combinatorics}.
Cambridge University Press, 2009.

\bibitem[Fran--70]{Fran--70}
J.\ Franks,
\emph{Anosov diffeomorphisms}.
In ``Global analysis (Proc.\ Sympos.\ Pure Math., Vol.\ XIV, Berkeley, Calif., 1968), pp. 61--93,
Amer.\ Math.\ Soc.\, 1970.

\bibitem[Gauss]{Gauss}
C.F.\ Gauss,
\emph{Werke, Vol.\ 2}.
G\"ottingen, 1876.


\bibitem[Grig--83]{Grig--83}
R.\ Grigorchuk,
\emph{On the Milnor problem of group growth}.
Dokl.\ Akad.\ Nauk SSSR \textbf{271} (1983), no.\ 1, 30--33.

\bibitem[Grig--84]{Grig--84}
R.\ Grigorchuk,
\emph{Degrees of growth of finitely generated groups 
and the theory of invariant means},
Izv.\ Akad.\ Nauk SSSR Ser.\ Mat.\ \textbf{48}(5) (1984), 939--985.
[In English Math.\ USSR Izv.\ \textbf{85} (1985) 259--300.]

\bibitem[Grig--14]{Grig--14}
R.\ Grigorchuk,
\emph{Milnor's problem on the growth of groups and its consequences}.
In ``Frontier in complex dynamics.
In celebration of John Milnor's 80th birthday.''
Princeton Mathematical Series \textbf{51}
(Princeton Univ.\ Press, 2014), 705--773.

\bibitem[GrHa--97]{GrHa--97}
R.\ Grigorchuk and P.\ de la Harpe,
\emph{On problems related to growth, entropy, and spectrum in group theory}.
J.\ Dynam.\ Control Systems \textbf{3} (1997), no.\ 1, 51--89.

\bibitem[Grom--81]{Grom--81} 
M.\ Gromov,
\emph{Groups of polynomial growth and expanding maps}.
Inst.\ Hautes \'Etudes Sci.\ Publ.\ Math.\ {\bf 53} (1981), 53--73.


\bibitem[Grom--87]{Grom--87} 
M.\ Gromov,
\emph{Hyperbolic groups}.
In ``Essays in group theory'', S.\ Gersten Editor,
Math.\ Sci.\ Res.\ Inst.\ Publ.\, \textbf{8} (Springer, 1987), 75--263.

\bibitem[Grom--93]{Grom--93} 
M.\ Gromov,
\emph{Asymptotic invariants of infinite groups}.
In ``Geometric group theory, Volume 2'',
London Math.\ Soc.\ Lecture Note Ser.\ {\bf 182}, Cambridge Univ.\ Press, 1993.


\bibitem[Guiv--73]{Guiv--73}
Y.\ Guivarc'h,
\emph{Croissance polynomiale et p\'eriodes des fonctions harmoniques}.
Bull.\ Soc.\ Math.\ France \textbf{101} (1973), 353--379.

\bibitem[GuRa--12]{GuRa--12}
Y.\ Guivarc'h and C.R.E.\ Raja,
\emph{Recurrence and ergodicity of random walks on linear groups
and homogeneous spaces}.
Ergodic Theory Danym.\ Systems \textbf{32} (2012), no.\ 4, 1313--1349.

\bibitem[Haus--14]{Haus--14}
F.\ Hausdorff,
\emph{Grundz\"uge der Mengenlehre}.
Veit \& Comp., 1914.

\bibitem[Hilb--90]{Hilb--90}
D.\ Hilbert,
\emph{Ueber die Theorie der algebraische Formen}.
Math.\ Ann.\ \textbf{36} (1890), no.~4, 473--534
[Gesammelte Abhandlungen, zweiter Band, 199--257].
See also \cite[Pages 143--224]{Hilbert}.

\bibitem[Hilbert]{Hilbert}
D.\ Hilbert,
\emph{Hilbert's invariant theory papers}.
Translated from the German by M.\ Ackerman.
With comments by R.\ Hermann.
Math Sci Press, Brookline, Mass., 1978.

\bibitem[Hjor--05]{Hjor--05}
G.\ Hjorth,
\emph{A converse to Dye's theorem}.
Trans.\ Amer.\ Math.\ Soc.\ \textbf{357} (2005),
3083--3103.


\bibitem[Jenk--73]{Jenk--73}
J.W.\ Jenkins,
\emph{A characterization of growth in locally compact groups}.
Bull.\ Amer.\ Math.\ Soc.\ \textbf{79} (1973), 103--106.

\bibitem[Kest--67]{Kest--67}
H.\ Kesten,
\emph{The Martin boundary of recurrent random walks on countable groups}.
Proc.\ Fifth Berkeley Sympos.\ M<th.\ Statist.\ and Probability
(Berkeley, Calif., 1965/66) Vol.\ II:
Contributions to Probability Theory, Part 2, pp.\ 51--74,
Univ.\ California, Berkeley, 1967.

\bibitem[Knut--97]{Knut--97}
D.E.\ Knuth,
\emph{The art of computer programming. Vol.\ 1}.
Third edition.
Addison--Wesley, 1997.

\bibitem[Krau--53]{Krau--53}
H.U.\ Krause,
\emph{Gruppenstruktur und Gruppenbild}.
Thesis, Eidgen\"ossische Technishe Hochschule, Z\`urich, 1953.
\tiny
\url{https://www.research-collection.ethz.ch/bitstream/handle/20.500.11850/133168/eth-32967-02.pdf?sequence=2&isAllowed=y}
\footnotesize

\bibitem[Lose--87]{Lose--87}
V.\ Losert,
\emph{On the structure of groups of polynomial growth}.
Math.] Z.\ \textbf{195} (1987), 109--117.

\bibitem[Mann--12]{Mann--12}
A.\ Mann,
\emph{How groups grow}.
London Math.\ Soc.\ Lecture Note Ser.\ {\bf 395},
Cambridge University Press, 2012.

\bibitem[Marg--67]{Marg--67}
G.A.\ Margulis,
\emph{Y-flows and three-dimensional manifolds}
(Appendix to \cite{AnSi--67}).
Russian Math.\ Surveys \textbf{22} (1967), no.\ 5, 164--166.

\bibitem[Marg--69]{Marg--69}
G.A.\ Margulis,
\emph{Applications of ergodic theory to the investigation of manifolds of negative curvature}.
Funct.\ Anal.\ Appl.\ \textbf{3} (1969), 335--336.

\bibitem[Miln--68a]{Miln--68a} 
J.\ Milnor,
\emph{Problem 5603}.
Amer.\ Math.\ Monthly \textbf{75} (1968), no.\ 6, 685.

\bibitem[Miln--68b]{Miln--68b} 
J.\ Milnor,
\emph{A note on curvature and fundamental group}.
J.\ Diff.\ Geom.\ \textbf{2} (1968), 1--7.
[Collected Papers, Volume I, 53 and 55--61.]

\bibitem[Miln--68c]{Miln--68c} 
J.\ Milnor,
\emph{Growth of finitely generated solvable groups}.
J.\ Diff.\ Geom.\ \textbf{2} (1968), 447--449.
[Collected Papers, Volume V, 155--157.]


\bibitem[vNeu--29]{vNeu--29}
J.\ von Neumann,
\emph{Zur allgemeinen Theorie des Masses}.
Fund.\ Math.\ \textbf{13} (1929) 73--116, 333 
[Collected works, Vol.\ I, 599--643].

\bibitem[OEIS]{OEIS}
The on-line encyclopedia of integer sequences,
\url{https://oeis.org}.

\bibitem[OrWe--80]{OrWe--80}
D.S.\ Ornstein and B.\ Weiss,
\emph{Ergodic theory of amenable group actions. I. The Rohlin lemma}.
Bull.\ Amer.\ Math.\ Soc.\ \textbf{2} (1980), 161--164.


\bibitem[Parr--92]{Parr--92}
W.\ Parry,
\emph{Growth series of some wreath products}.
Trans.\ Amer.\ Math.\ Soc.\ \textbf{331} (1992), no.\ 2, 751--759.

\bibitem[PlTh--72]{PlTh--72}
J.F.\ Plante and W.P.\ Thurston,
\emph{Anosov flows and the fundamental group}.
Topology \textbf{11} (1972), 147--150.

\bibitem[Poin--82]{Poin--82}
H.\ Poincar\'e,
\emph{Th\'eorie des groupes fuchsiens},
Acta Math.\ \textbf{1} (1982), 1--62.
[Oeuvres, tome II, pages 108--168.]

\bibitem[Poly--21]{Poly--21}
G.\ P\'olya,
\emph{\"Uber eine Aufgabe der Wahrscheinlichkeitsrechnung betreffend
die Irr\-fahrt im Strassennetz}.
Math.\ Ann.\ \textbf{84} (1921), no.\ 1--2, 149--160.

\bibitem[Rose--74]{Rose--74}
J.M.\ Rosenblatt,
\emph{Invariant measures and growth conditions}.
Trans.\ Amer.\ Math.\ Soc.\ \textbf{193} (1974), 33--53. 

\bibitem[Serr--71]{Serr--71}
J-P.\ Serre,
\emph{Cohomologie des groupes discrets}.
In ``Prospects in mathematics (Princeton, 1970)'', 
Annals of Math.\ Studies \textbf{70} (Princeton Univ.\ Press, 1971), 77--169.

\bibitem[Sier--46]{Sier--46}
W.\ Sierpinski,
\emph{Sur la non-existence des d\'ecompositions paradoxales d'ensembles lin\'eaires}.
Actas Acad.\ Ci.\ Lima \textbf{9} (1946), 113--117.

\bibitem[Sylvester]{Sylvester}
\emph{The collected mathematical papers of James Joseph Sylvester},
Volume III (1870--1883).
Cambridge University Press, 1909.

\bibitem[Stan--96]{Stan--96}
R.\ Stanley,
\emph{Combinatorics and commutative algebra}, Second edition.
Birkh\"auser, 1996
[First edition 1983].

\bibitem[Stan--15]{Stan--15}
R.\ Stanley,
\emph{Catalan numbers},
with an appendix on \emph{History of Catalan numbers} by I.\ Pak.
Cambridge University Press, 2015.
[See also R.\ Stanley, \emph{Catalan addendum}, on
\url{http://www-math.mit.edu/~rstan/ec/catadd.pdf}
and I.\ Pak, \emph{Catalan numbers}, on
\url{https://www.math.ucla.edu/~pak/lectures/Cat/pakcat.htm}].

\bibitem[Stol--96]{Stol--96}
M.\ Stoll,
\emph{Rational and transcendental growth series for the higher Heisenberg groups}.
Invent.\ Math.\ \textbf{126} (1996), no.\ 1, 85--109.


\bibitem[Svar--55]{Svar--55}
A.S.\ \v{S}varc,
\emph{A volume invariant of coverings}.
Dokl.\ Akad.\ Nauk.\ SSSR \textbf{105} (1955), 32--34.

\bibitem[Svar--08]{Svar--08}
A. Schwarz [= A.S.\ \v{S}varc],
\emph{My life in science}.
Updated in November 2008,
in Albert Schwarz's Homepage -- UC Davis Math,
\url{https://web.archive.org/web/20101125124406/http://www.math.ucdavis.edu/~schwarz/}

\bibitem[Tars--36]{Tars--36}
A.\ Tarski,
\emph{Algebraische Fassung des Massproblems}.
Fund.\ Math.\ \textbf{31} (1938), 47--66
[Collected Papers, Volume 2, 453--472].

\bibitem[VaSC--92]{VaSC--92}
N.Th.\ Varopoulos, L.\ Saloff-Coste, and T.\ Coulhon,
\emph{Analysis and geometry on groups}.
Cambridge Univ.\ Press, 1992.

\bibitem[Wagr--82]{Wagr--82}
P.\ Wagreich.
\emph{The growth function of a discrete group}.
In ``Group actions and vector fields (Vancouver, B.C, 1982)'',
Lecture Notes in Math.\ \textbf{956} (Springer, 1982), 125--144.

\bibitem[Wolf--68]{Wolf--68}
J.A.\ Wolf,
\emph{Growth of finitely generated solvable groups and curvature
of Riemannian manifolds}.
J.\ Diff.\ Geom.\ \textbf{2} (1968), 421--446.

\end{thebibliography}
\end{document}